\documentclass[12pt]{article}

\usepackage{amssymb}
\usepackage{amsmath}
\usepackage[mathscr]{eucal}
\usepackage{cite}

\begin{document}

\renewcommand{\citeleft}{{\rm [}}
\renewcommand{\citeright}{{\rm ]}}
\renewcommand{\citepunct}{{\rm,\ }}
\renewcommand{\citemid}{{\rm,\ }}

\newcounter{abschnitt}
\newtheorem{theorem}{Theorem}[abschnitt]
\newtheorem{koro}{Corollary}[abschnitt]
\newtheorem{defi}{Definition}
\newtheorem{satz}[koro]{Theorem}
\newtheorem{lem}[koro]{Lemma}

\newcounter{saveeqn}
\newcommand{\alpheqn}{\setcounter{saveeqn}{\value{abschnitt}}
\renewcommand{\theequation}{\mbox{\arabic{saveeqn}.\arabic{equation}}}}
\newcommand{\reseteqn}{\setcounter{equation}{0}
\renewcommand{\theequation}{\arabic{equation}}}

\hyphenation{convex} \hyphenation{bodies}

\sloppy

\phantom{a}

\vspace{-1.5cm}

\begin{center} \begin{LARGE} $\mathrm{GL}(n)$ contravariant Minkowski
valuations \\[0.7cm] \end{LARGE}

\begin{large} Franz E. Schuster and Thomas Wannerer \end{large}
\end{center}

\vspace{-0.8cm}

\begin{quote}
\footnotesize{ \vskip 1truecm\noindent {\bf Abstract.} A complete
classification of all continuous $\mathrm{GL}(n)$ contravariant
Minkowski valuations is established. As an application we present
a family of sharp isoperimetric inequalities for such valuations
which generalize the classical Petty projection inequality.}
\end{quote}

\vspace{0.6cm}

\centerline{\large{\bf{ \setcounter{abschnitt}{1}
\arabic{abschnitt}. Introduction}}} \alpheqn

\vspace{0.6cm}

The projection body of a convex body is one of the central
notions that Minkowski introduced within convex geometry. Over the
past four decades it has become evident that the projection
operator, its range (i.e.\ the class of centered zonoids), and
its polar are objects which arise naturally in a number of
different areas, see e.g., \textbf{\cite{gardner2ed, ludwig02,
LYZ2004, LYZ2010, Schneider:Weil83, schnweil, steineder08,
Thompson}}. The most important affine isoperimetric inequality
for projection bodies is the {\it Petty projection inequality}
\textbf{\cite{petty67}} which is significantly stronger than the
classical isoperimetric inequality. This remarkable inequality
also forms the geometric core of the affine Zhang--Sobolev
inequality \textbf{\cite{LYZ2002b, Zhang99}} which strengthens
and directly implies the classical Sobolev inequality.

The special role of projection bodies in affine convex geometry
was \linebreak demonstrated only recently by two results of Ludwig
\textbf{\cite{ludwig02, Ludwig:Minkowski}}: The projection
\linebreak operator was characterized as the unique continuous
Minkowski valuation which is $\mathrm{GL}(n)$ contravariant {\it
and} invariant under translations. Moreover, the \linebreak
assumption of translation invariance can be omitted when the
domain of the valuation is restricted to convex bodies containing
the origin. Through the seminal work of Ludwig, convex and star
body valued valuations have become \linebreak the focus of
increased attention, see e.g., \textbf{\cite{haberl08, hablud06,
haberl09, Ludwig06, Schu09}}. A very recent \linebreak development
in this area explores the connections between these valuations and
the theory of isoperimetric inequalities, see e.g.,
\textbf{\cite{habschu09, hs2}}.

In this article we establish a complete classification of all
continuous and $\mathrm{GL}(n)$ contravariant Minkowski
valuations, without any further assumption on their behavior under
translations or any restrictions on their domain. We show that
there is a two parameter family of such valuations generated by
the projection operator and the convex hull with the origin. We
also obtain a Petty projection type inequality for each member of
this family and identify {\it the} Petty projection inequality as
the strongest inequality.

\pagebreak

Let $\mathcal{K}^n$ denote the space of convex bodies (compact
convex sets) in $\mathbb{R}^n$, $n \geq 3$, endowed with the
Hausdorff metric and let $\mathcal{K}_0^n$ denote the subset of
$\mathcal{K}^n$ of bodies containing the origin. A convex body $K$
is uniquely determined by its support function
$h(K,x)=\max\{x\cdot y: y \in K\}$, for $x \in \mathbb{R}^n$.

\vspace{0.3cm}

\noindent {\bf Definition} \emph{A map $\Phi: \mathcal{K}^n
\rightarrow \mathcal{K}^n$ is called a Minkowski valuation if
\[\Phi K + \Phi L = \Phi(K \cup L) + \Phi(K \cap L),  \]
whenever $K, L, K \cup L \in \mathcal{K}^n$ and addition on
$\mathcal{K}^n$ is Minkowski addition. \\[0.1cm]
The map $\Phi$ is called $\mathrm{GL}(n)$ contravariant if there
exists  $q \in \mathbb{R}$ such that for all $\phi \in
\mathrm{GL}(n)$ and all $K \in \mathcal{K}^n$,
\[\Phi(\phi K) = |\det \phi|^q \phi^{-\mathrm{T}}\Phi K,  \]
where $\phi^{-\mathrm{T}}$ denotes the transpose of the inverse
of $\phi$.}

\vspace{0.3cm}

The classical example of a $\mathrm{GL}(n)$ contravariant
Minkowski valuation is the projection operator: The {\it
projection body} $\Pi K$ of $K \in \mathcal{K}^n$ is defined by
\begin{equation} \label{defpi}
h(\Pi K,u)=\mathrm{vol}_{n-1}(K|u^{\bot}), \qquad u \in S^{n-1},
\end{equation}
where $K|u^{\bot}$ denotes the projection of $K$ onto the
hyperplane orthogonal to $u$. \linebreak It was first proved by
Petty \textbf{\cite{petty67b}} that for all $\phi \in
\mathrm{GL}(n)$ and all $K \in \mathcal{K}^n$,
\[\Pi(\phi K) = |\det \phi | \phi^{-\mathrm{T}}\, \Pi K.  \]

The notion of valuation plays a central role in geometry. It was
the critical ingredient in the solution of Hilbert's Third Problem
and has since been intimately tied to the dissection theory of
polytopes (see \textbf{\cite{McMullen93}}). Over the last decades
the theory of valuations has evolved enormously and had a
tremendous impact on various disciplines, see e.g.,
\textbf{\cite{Alesker99, Alesker01, Bernig07, Bernig03,
bernigfu10, fu06, klain99, klain00, Klain:Rota, Ludwig:matrix,
centro}}.

First results on Minkowski valuations which are {\it rigid motion}
equivariant were obtained by Schneider \textbf{\cite{schneider74}}
in the 1970s (see \textbf{\cite{kiderlen05, schnschu, Schu06a,
Schu09}} for recent extensions of these results). The starting
point for the systematic study of convex and star body valued
valuations which are compatible with {\it linear transformations}
were two highly influential articles by Ludwig
\textbf{\cite{ludwig02, Ludwig:Minkowski}}. The following theorem
is an example of the numerous results obtained there:

\vspace{0.3cm}

\noindent {\bf Theorem 1 (Ludwig \cite{Ludwig:Minkowski})} \emph{A
map $\Phi :\mathcal{K}_0^n \rightarrow \mathcal{K}^n$ is a
$\mathrm{GL}(n)$ contravariant continuous Minkowski valuation if
and only if there exists a constant $c \geq 0$ such that for
every $K\in\mathcal{K}^n_0$,}
\[\Phi K=c\,\Pi K.\]

\vspace{0.3cm}

If we consider Minkowski valuations defined on all convex bodies,
then the projection operator is no longer the only
$\mathrm{GL}(n)$ contravariant valuation. To see this consider
the map $\Pi_{\mathrm{o}}:\mathcal{K}^n \rightarrow \mathcal{K}^n$
defined by
\[\Pi_{\mathrm{o}}K := \Pi(\mathrm{conv}(\{0\}\cup K)),\]
where $\mathrm{conv}(\{0\}\cup K)$ denotes the convex hull of $K$
and the origin. It is easy to see that $\Pi_{\mathrm{o}}$ is a
$\mathrm{GL}(n)$ contravariant continuous Minkowski valuation
but, clearly, it is not a multiple of the projection operator.

The main object of this paper is to show that all
$\mathrm{GL}(n)$ contravariant continuous Minkowski valuations on
$\mathcal{K}^n$ are given by combinations of the projection
operator and the map $\Pi_{\mathrm{o}}$.

\vspace{0.4cm}

\noindent {\bf Theorem \!2} \emph{A map $\Phi :\mathcal{K}^n
\rightarrow \mathcal{K}^n$ is a $\mathrm{GL}(n)$ contravariant
continuous Minkowski valuation if and only if there exist
constants $c_1, c_2 \geq 0$ such that for every $K \in
\mathcal{K}^n$,}
\[\Phi K=c_1\:\!\Pi K + c_2\;\!\Pi_{\mathrm{o}}K.\]

\vspace{0.4cm}

Let $K^*=\{x \in \mathbb{R}^n: x \cdot y \leq 1 \mbox{ for all } y
\in K\}$ denote the {\it polar body} of a convex body $K$
containing the origin in its interior. We write $\Phi^* K$ to
denote $(\Phi K)^*$, $V(K)$ for the volume of $K$, and $B$ for
the Euclidean unit ball.

In the early 1970s Petty established a fundamental affine
isoperimetric inequality, now known as the Petty projection
inequality. (For important recent generalizations, see
\textbf{\cite{LYZ2000a, LYZ2010}}.) If $K \in \mathcal{K}^n$ has
nonempty interior, then
\begin{equation} \label{petty1}
V(K)^{n-1}V(\Pi^*K) \leq V(B)^{n-1}V(\Pi^*B)
\end{equation}
with equality if and only if $K$ is an ellipsoid.

As an application of Theorem 2 we extend the Petty projection
inequality to the entire class of $\mathrm{GL}(n)$ contravariant
continuous Minkowski valuations which are {\it nontrivial}, i.e.,
which do not map every convex body to the origin.

\vspace{0.4cm}

\noindent {\bf Theorem 3} \emph{Let $K \in \mathcal{K}^n$ have
nonempty interior. If $\Phi :\mathcal{K}^n \rightarrow
\mathcal{K}^n$ is a nontrivial $\mathrm{GL}(n)$ contravariant
continuous Minkowski valuation, then
\begin{equation*}
V(K)^{n-1}V(\Phi^*K) \leq V(B)^{n-1}V(\Phi^*B).
\end{equation*}
If $\Phi$ is not a multiple of $\Pi$, there is equality if and
only if $K$ is an ellipsoid containing the origin; otherwise
equality holds if and only if $K$ is an ellipsoid.}

\vspace{0.4cm}

In view of Theorems 2 and 3, the natural problem arises to
determine for fixed $K \in \mathcal{K}^n$ the maximum value of
$V(\Phi^*K)$ among all suitably normalized (say, e.g., $\Phi B =
\Pi B$) $\mathrm{GL}(n)$ contravariant continuous Minkowski
valuations. Here, we will show that
\[V(\Phi^*K) \leq V(\Pi^*K).   \]
This shows that the classical Petty projection inequality gives
rise to the strongest inequality among the inequalities of
Theorem 3.

\vspace{1.2cm}

\centerline{\large{\bf{ \setcounter{abschnitt}{2}
\arabic{abschnitt}. Background material}}}

\vspace{0.7cm} \reseteqn \alpheqn

In the following we state for quick reference some basic facts
about convex bodies and the geometric inequalities needed in the
proof of Theorem 3. For general reference the reader may wish to
consult the book by Schneider \textbf{\cite{schneider93}}.

A convex body $K \in \mathcal{K}^n$ is uniquely determined by the
values of its \linebreak support function $h(K,\cdot)$ on
$S^{n-1}$. If $\phi \in \mathrm{GL}(n)$ and $K \in
\mathcal{K}^n$, then for every $x \in \mathbb{R}^n$,
\begin{equation} \label{supptrans}
h(\phi K,x) = h(K,\phi^{\,\mathrm{T}}x).
\end{equation}

A convex body $K \in \mathcal{K}^n$ with nonempty interior is
also determined up to translations by its {\it surface area
measure} $S_{n-1}(K,\cdot)$. Recall that for a Borel set $\omega
\subseteq S^{n-1}$, $S_{n-1}(K,\omega)$ is the $(n-1)$-dimensional
Hausdorff measure of the set of all boundary points of $K$ at
which there exists a normal vector of $K$ belonging to $\omega$.

It is a well known fact that the projection function
$\mathrm{vol}_{n-1}(K|\,\cdot^{\bot})$ of a convex body $K \in
\mathcal{K}^n$ is (up to a constant) the {\it cosine transform}
of the surface area measure of $K$. More precisely, we have
\begin{equation} \label{projform}
\mathrm{vol}_{n-1}(K|u^{\bot}) =\frac{1}{2} \int_{S^{n-1}}|u \cdot v|\,dS_{n-1}(K,v), \qquad u \in S^{n-1}.
\end{equation}

For a compact set $L$ in $\mathbb{R}^n$ which is star-shaped with
respect to the origin, its radial function is defined by
$\rho(L,x)=\max\{\lambda \geq 0: \lambda x \in L\}$, $x \in
\mathbb{R}^n\backslash \{0\}$. For $\alpha, \beta \geq 0$ (not
both zero) and $K, L \in \mathcal{K}^n$ containing the origin in
their interiors, the harmonic (radial) combination $\alpha \cdot
K \, \widehat{+} \, \beta \cdot L$ of $K$ and $L$ is the convex
body defined by
\begin{equation*}
\rho(\alpha \cdot K \, \widehat{+} \, \beta \cdot L,\cdot)^{-1} =
\alpha \rho(K,\cdot)^{-1} + \beta \rho(L,\cdot)^{-1}.
\end{equation*}

If $K \in \mathcal{K}^n$ is a convex body containing the origin
in its interior, then it follows from the definitions of support
and radial functions, and the definition of the polar body of
$K$, that
\begin{equation*}
\rho(K^*,\cdot)=h(K,\cdot)^{-1} \qquad \mbox{and} \qquad h(K^*,\cdot)=\rho(K,\cdot)^{-1}.
\end{equation*}
Thus, we have
\begin{equation} \label{harsum}
\alpha \cdot K \, \widehat{+} \, \beta \cdot L = (\alpha K^* +
\beta L^*)^*,
\end{equation}
which shows that the set of convex bodies containing the origin
in their interiors is closed under harmonic combinations.

First results on harmonic combinations of convex bodies were
obtained in the early 1960's by Firey \textbf{\cite{firey61}}. In
particular, he established the following {\it harmonic dual
Brunn--Minkowski inequality}: If $K,L \in \mathcal{K}^n$ contain
the origin in their interiors, then
\begin{equation} \label{dualbm}
V(K \, \widehat{+} \, L)^{-1/n} \geq V(K)^{-1/n} + V(L)^{-1/n},
\end{equation}
with equality if and only if $K$ and $L$ are dilates.

A map $\Phi: \mathcal{K}^n \rightarrow \mathcal{K}^n$ is called
$\mathrm{SL}(n)$ {\it contravariant} if for all $\phi \in
\mathrm{SL}(n)$ and all $K \in \mathcal{K}^n$,
\[\Phi(\phi K) = \phi^{-\mathrm{T}}\Phi K,  \]
and it is called (positively) {\it homogeneous} of degree $r$, $r
\in \mathbb{R}$, if for all $\lambda
> 0$ and all $K \in \mathcal{K}^n$,
\[\Phi(\lambda K)=\lambda^r \Phi K.  \]

If $\Phi$ is $\mathrm{SL}(n)$ contravariant and homogeneous of
degree $r$, then we have for all $K \in \mathcal{K}^n$,
\begin{equation} \label{eq:1.2}
\Phi(\phi K) = (\det \phi)^{(r+1)/n} \phi^{-\mathrm{T}}\Phi K,
\end{equation}
for every $\phi \in \mathrm{GL}(n)$ with $\det \phi > 0$.

Since, clearly, every $\mathrm{GL}(n)$ contravariant map is
$\mathrm{SL}(n)$ contravariant and homogeneous, the following
result is a stronger version of Theorem 1.

\vspace{0.3cm}

\noindent {\bf Theorem 4 (Ludwig \cite{Ludwig:Minkowski})} \emph{A
map $\Phi :\mathcal{K}_0^n \rightarrow \mathcal{K}^n$ is an
$\mathrm{SL}(n)$ contravariant and homogeneous continuous
Minkowski valuation if and only if there exists a constant $c
\geq 0$ such that for every $K\in\mathcal{K}^n_0$,}
\[\Phi K=c\,\Pi K.\]

\vspace{1cm}

\centerline{\large{\bf{ \setcounter{abschnitt}{3}
\arabic{abschnitt}. An auxiliary result}}}

\vspace{0.7cm} \reseteqn \alpheqn \setcounter{koro}{0}

In this section we show that any $\mathrm{SL}(n)$ contravariant
continuous Minkowski valuation which is homogeneous of degree $r
\neq n - 1$ is trivial. \linebreak The proof is based on the
ideas and techniques developed by Ludwig
\textbf{\cite{Ludwig:Minkowski}}.

As a first step we consider the image of convex bodies contained
in a hyperplane. We denote by $\mathrm{aff}\,K$ the affine hull
of $K \subseteq \mathbb{R}^n$.

\vspace{0.4cm}

\noindent {\bf Lemma 5} \emph{ Let
$\Phi:\mathcal{K}^n\rightarrow\mathcal{K}^n$ be $\mathrm{SL}(n)$
contravariant and homogeneous of degree $r$. Then the following
statements hold:
\begin{itemize}
\item[(i)] If $\dim K < n-2$,  then $\Phi K=\{0\}$.
\item[(ii)] If $\dim K = n-2$, then
\[\left \{ \begin{array}{ll} \Phi K = \{0\} & \mbox{if }  0 \in \mathrm{aff}\,
K,
\\ \Phi K \subseteq \mathrm{aff}(\{0\}\cup K)^\bot & \mbox{if }
0 \notin \mathrm{aff}\, K. \end{array} \right .   \]
\item[(iii)] If $\dim K = n - 1$ and $0 \in \mathrm{aff}\, K$, then
$\Phi K\subseteq (\mathrm{aff} K)^\perp$.
    \item[(iv)]
If  $K$ is contained in a hyperplane through the origin and $r\neq
n-1$, then $\Phi K=\{0\}$.
\end{itemize}
}

\noindent {\it Proof}\,: Let $\{e_1, \ldots, e_n\}$ be a fixed
orthonormal basis of $\mathbb{R}^n$. For $k \geq 1$, we use the
direct sum decomposition $\mathbb{R}^n = \mathbb{R}^k \oplus
\mathbb{R}^{n-k}$, where
\[\mathbb{R}^k = \mathrm{span}\,\{e_1, \ldots, e_k\} \qquad \mbox{and} \qquad \mathbb{R}^{n-k} = \mathrm{span}\,\{e_{k+1}, \ldots, e_n\}.  \]
Let $\phi \in \mathrm{SL}(n)$ be the matrix defined (with respect
to the chosen basis) by
\[\phi = \begin{pmatrix}I_k & B \\ 0 & A\end{pmatrix},\]
where $I_k$ is the $k\times k$ identity matrix, $B$ is an
arbitrary $k \times (n-k)$ matrix, and $A \in \mathrm{SL}(n-k)$.
A simple calculation shows that
\begin{equation} \label{phitrans}
\phi^{-\mathrm{T}}=\begin{pmatrix} I_k & 0\\
-A^{-\mathrm{T}}B^{\mathrm{T}} & A^{-\mathrm{T}}\end{pmatrix}.
\end{equation}
If $K \in \mathcal{K}^n$ and $K \subseteq \mathbb{R}^k$, we have
\begin{equation}\label{eq:3.1}
\phi K = K.
\end{equation}
Now let $x \in \Phi K$ and write $x = x' + x''$, where $x' \in
\mathbb{R}^k$, $x'' \in \mathbb{R}^{n-k}$. Since $\Phi$ is
$\mathrm{SL}(n)$ contravariant, it follows from (\ref{phitrans})
and (\ref{eq:3.1}) that
\begin{equation}\label{eq:3.2}
\phi^{-\mathrm{T}}x =  x' -A^{-\mathrm{T}}(B^{\mathrm{T}} x'- x'')
\in \Phi K.
\end{equation}
This holds for every $A \in \mathrm{SL}(n-k)$ and every $k \times
(n - k)$ matrix $B$. Since $\Phi K$ is bounded, this implies that
$x' = 0$. Thus, we have
\begin{equation} \label{cont17}
\Phi K \subseteq \mathbb{R}^{n-k}.
\end{equation}
Moreover, if $n - k \geq 2$, it follows from (\ref{eq:3.2}) that
also $x''=0$, i.e., $\Phi K = \{0\}$.

Since every convex body $K \in \mathcal{K}^n$ with $\dim K < n -
2$ or $\dim K = n - 2$ and $0 \in \mathrm{aff}\,K$ is the linear
image of some convex body contained in $\mathbb{R}^{n-2}$, we
obtain (i) and the first assertion of (ii). The second assertion
of (ii) and statement (iii) are immediate consequences of
(\ref{cont17}).

Finally, let $K \subseteq \mathbb{R}^{n-1}$ and $r \neq n-1$. If
$\psi \in \mathrm{GL}(n)$ is defined by
\[\psi=\begin{pmatrix} I_{n-1} & 0 \\ 0 & s \end{pmatrix},\]
where $s > 0$, then, by (\ref{eq:1.2}),
\[\Phi K = \Phi (\psi K) = s^{(r-(n-1))/n} \Phi K.\]
Since this holds for every $s > 0$ and $\Phi K$ is bounded, we
deduce (iv). \hfill $\blacksquare$

\vspace{0.4cm}

Our next result reduces the proof of Theorem 2 to Minkowski
valuations which are homogeneous of degree $n - 1$.

\vspace{0.4cm}

\noindent {\bf Proposition 6} \emph{If $\Phi: \mathcal{K}^n
\rightarrow \mathcal{K}^n$ is an $\mathrm{SL}(n)$ contravariant
continuous Minkowski valuation which is homogeneous of degree
$r\neq n-1$, then for every $K\in\mathcal{K}^n$,}
\[\Phi K = \{0\}.\]

\vspace{0.2cm}

\noindent {\it Proof}\,: Let $\{e_1, \ldots, e_n\}$ be a fixed
orthonormal basis of $\mathbb{R}^n$. For $0<\lambda<1$ and $1\leq
i<j\leq n$, we denote by $H_\lambda=H_\lambda(i,j)$ the
hyperplane through the origin with normal vector $\lambda e_j -
(1-\lambda)e_i $, that is,
\[H_\lambda=\mathrm{span}(\{\lambda e_i+(1-\lambda)e_j\}\cup\{e_k:k\neq i,j\}).\]
Let $S$ be the $(n-1)$-dimensional simplex with vertices
$\{e_1,\ldots,e_n\}$. Then $H_\lambda$ dissects $S$ into two
simplices $S\cap H^+_\lambda$ and $S\cap H^-_\lambda$, where
$H^+_\lambda$, $H^-_\lambda$ denote the closed halfspaces bounded
by $H_\lambda$.  More precisely, we have
\begin{eqnarray*}
S \cap H^+_\lambda & = & \mathrm{conv}(\{\lambda
e_i+(1-\lambda)e_j\}\cup\{e_k:k\neq i\}), \\
S \cap H^-_\lambda & = & \mathrm{conv}(\{\lambda
e_i+(1-\lambda)e_j\}\cup\{e_k:k\neq j\}).
\end{eqnarray*}
Define linear maps $\phi_\lambda=\phi_\lambda(i,j)$ and
$\psi_\lambda=\phi_\lambda(i,j)$ by
$$\phi_\lambda e_i= \lambda e_i +(1-\lambda)e_j, \ \ \phi_\lambda e_k= e_k\ \text{ for } k\neq i,$$
$$\psi_\lambda e_j= \lambda e_i +(1-\lambda)e_j, \ \ \phi_\lambda e_k= e_k\ \text{ for } k\neq j.$$
Then it is easy to see that
\begin{equation} \label{linhalf17}
S \cap H^+_\lambda=\phi_\lambda S \qquad \mbox{and} \qquad S\cap
H^-_\lambda=\psi_\lambda S.
\end{equation}

Since $\Phi$ is a Minkowski valuation, we have
\begin{equation} \label{valprop17}
\Phi S + \Phi(S \cap H_\lambda) = \Phi(S \cap H^+_\lambda) +
\Phi(S\cap H^-_\lambda).
\end{equation}
Since $\Phi$ is also $\mathrm{SL}(n)$ contravariant and
homogeneous of degree $r \neq n - 1$, it follows from Lemma 5
(iv) that $\Phi(S \cap H_{\lambda}) = \{0\}$. Therefore,
(\ref{eq:1.2}), (\ref{linhalf17}), and (\ref{valprop17}) yield
\begin{equation} \label{eq:4.2}
\Phi S = \lambda^q \phi^{-\mathrm{T}}_\lambda \Phi S+
(1-\lambda)^q \psi^{-\mathrm{T}}_\lambda \Phi S,
\end{equation}
where $q=(r+1)/n$.

Now let $1\leq k \leq n$ and choose $1 \leq i < j \leq n$ such
that $k \neq i,j$. This is possible since $n \geq 3$. By
(\ref{supptrans}) and (\ref{eq:4.2}), we have
\[h(\Phi S,e_k)=\lambda^q h(\Phi S,e_k)+(1-\lambda)^q  h(\Phi
S,e_k),\] for every $0 < \lambda < 1$. Since $\Phi$ is
homogeneous of degree $r \neq n-1$, we have $q \neq 1$ which
implies $h(\Phi S,e_k)=0$. Similarly, we obtain $h(\Phi
S,-e_k)=0$. Since this holds for every $1\leq k \leq n$, we must
have
\begin{equation} \label{simnull}
\Phi S=\{0\}.
\end{equation}
Let $T$ be an arbitrary $(n-1)$-dimensional simplex. If $T$ is contained in a hyperplane through the origin, then, by Lemma 5 (iv), $\Phi\, T = \{0\}$. Otherwise there exists a linear transformation
$\phi \in \mathrm{GL}(n)$ with $\det \phi > 0$ such that $\phi S = T$. Thus, by (\ref{eq:1.2}) and (\ref{simnull}), we also obtain
$\Phi\, T = \{0\}$.

Now let $F \in \mathcal{K}^n$ be an $(n-1)$-dimensional polytope.
We dissect $F$ into $(n-1)$-dimensional simplices $S_i$, $i = 1,
\ldots, m$, that is, $F = S_1 \cup \ldots \cup S_m$ and $\dim S_i
\cap S_j \leq n - 2$ whenever $i \neq j$. For $x \in
\mathbb{R}^n$ fixed, we define a real valued valuation
$\varphi_x: \mathcal{K}^n \rightarrow \mathbb{R}$ by
\[\varphi_x(K)=h(\Phi K,x).  \]
From the continuity of $\Phi$, it follows that $\varphi_x$ is also continuous. Therefore, $\varphi_x$ satisfies
the inclusion-exclusion principle
\[\varphi_x(F)=\sum_I (-1)^{|I|-1}\varphi_x(S_I),    \]
where the sum is taken over all ordered $k$-tuples $I =
(i_1,\ldots,i_k)$ such that $1 \leq i_1 < \ldots < i_k \leq n$
and $k = 1, \ldots ,m$. Here $|I|$ denotes the cardinality of $I$
and $S_I = S_{i_1} \cap \ldots \cap S_{i_k}$ (cf.\
\textbf{\cite[\textnormal{p.\ 7}]{Klain:Rota}}). By Lemma 5 and
the fact that $\Phi S_i = \{0\}$ for $i = 1, \ldots, m$, it
follows that $\varphi_x(F)=0$ for every $x \in \mathbb{R}^n$.
Consequently,
\begin{equation} \label{phifnull}
\Phi F = \{0\}
\end{equation}
for every $(n-1)$-dimensional polytope $F$.

Next, let $P \in \mathcal{K}^n$ be an $n$-dimensional polytope.
If $0\in P$, then $\Phi P = \{0\}$ by Theorem 4. Therefore,
suppose that $0\notin P$. We call a facet $F$ of $P$
\emph{visible} if $[0,x)\cap P=\emptyset$ for every $x\in F$. If
$P$ has exactly one visible facet $F$, then, since $\Phi$ is a
Minkowski valuation, we have
\[\Phi P_{\mathrm{o}} + \Phi F = \Phi F_{\mathrm{o}} + \Phi P,\]
where we use $K_{\mathrm{o}}$ to denote $\mathrm{conv}(\{0\} \cup
K)$. Since, by Theorem 4, $\Phi K = \{0\}$ for all $K \in
\mathcal{K}^n_{\mathrm{o}}$ and $\Phi F=\{0\}$ by
(\ref{phifnull}), we obtain $\Phi P=\{0\}$. If $P$ has $m > 1$
visible facets $F_1,\ldots,F_m$, we define polytopes $C_i \in
\mathcal{K}^n$, $i = 1, \ldots, m$, by
\[C_i = P \cap \bigcup_{t\geq0} t F_i.  \]
Clearly, $P = C_1 \cup \ldots \cup C_m$ and $\dim C_i\cap C_j \leq n - 1$ for $i\neq j$. Moreover, each $C_i$, $i = 1, \ldots, m$, has exactly one visible facet.
Hence, $\Phi\, C_i = \{0\}$ for $i = 1, \ldots, m$. Thus, as before, the inclusion-exclusion principle and (\ref{phifnull}) imply
\begin{equation} \label{phipnull}
\Phi P = \{0\}
\end{equation}
for every $n$-dimensional polytope $P$.

Finally, a combination of Lemma 5, (\ref{phifnull}), and
(\ref{phipnull}) shows that the map $\Phi$ vanishes on all
polytopes. By the continuity of $\Phi$, the same holds true for
all
 convex bodies. \hfill $\blacksquare$

\pagebreak

\centerline{\large{\bf{ \setcounter{abschnitt}{4}
\arabic{abschnitt}. Proof of the main result}}}

\vspace{0.7cm} \reseteqn \alpheqn \setcounter{koro}{0}

After these preparations, we are now in a position to present the
proof of Theorem 2. In fact, since $\mathrm{GL}(n)$ contravariant
Minkowski valuations are $\mathrm{SL}(n)$ contravariant and
homogeneous, we prove a slightly stronger result, analogous to
Theorem 4.

\vspace{0.4cm}

\noindent {\bf Theorem 7} \emph{A map $\Phi\!\: :\mathcal{K}^n
\rightarrow \mathcal{K}^n$ is an $\mathrm{SL}(n)$ contravariant
and homogeneous continuous Minkowski valuation if and only if
there exist constants $c_1, c_2 \geq 0$ such that for every $K
\in \mathcal{K}^n$,}
\[\Phi K=c_1\:\!\Pi K + c_2\;\!\Pi_{\mathrm{o}}K.\]

\vspace{0.2cm}

\noindent {\it Proof}\,: By Proposition 6, it is sufficient to
consider $\Phi$ that are homogeneous of degree $n - 1$. We first
show that in this case, $\Phi$ is already determined by its values
on convex bodies of dimension $n-1$. To this end, let $\{e_1,
\ldots, e_n\}$ be an orthonormal basis of $\mathbb{R}^n$ and
define $\phi_s \in \mathrm{GL}(n)$ by
\[\phi_s=\begin{pmatrix} s\,I_{n-1} & 0 \\  0 & 1\end{pmatrix},\]
where $s > 0$ and, as before, $I_k$ denotes the $k \times k$
identity matrix. Since $\Phi$ is homogeneous of degree $n - 1$,
(\ref{eq:1.2}) yields
\[\phi_s \Phi K = \Phi (s\phi_s^{-\mathrm{T}} K).\]
Hence, by the continuity of $\Phi$, letting $s\rightarrow 0$, we have
\[\phi_0 \Phi K =\Phi(K|e_n^{\bot}).\]
Since $h(\phi_0 L,e_n) = h(L,e_n)$ for every $L \in \mathcal{K}^n$, we obtain
\begin{equation}\label{eq:3.3}h(\Phi K,e_n) = h(\Phi(K|e_n^{\bot}),e_n). \end{equation}
Now let $u\in S^{n-1}$ and choose $\vartheta \in \mathrm{SO}(n)$
such that $\vartheta e_n= u$. Using (\ref{supptrans}),
(\ref{eq:1.2}) and (\ref{eq:3.3}), we conclude that
\begin{equation} \label{detlow}
h(\Phi K,u) = h(\Phi((\vartheta^{-1} K)|e_n^{\bot}),e_n)=h(\Phi(K|u^{\bot}),u).
\end{equation}

\vspace{0.1cm}

Next, we show that  there exists a constant $c\geq0$ such that
\begin{equation}\label{eq:3.4}\Phi K=c\, \mathrm{vol}_{n-1}(K_{\mathrm{o}})\,[-u,u]\end{equation}
for every $K\in\mathcal{K}^n$ contained in $u^\perp$, $u\in
S^{n-1}$, with $\dim K\leq n-2$. Here, $[-u,u]$ denotes the
segment with endpoints $-u, u$ and
$K_{\mathrm{o}}=\mathrm{conv}\,(\{0\} \cup K)$.

\pagebreak

If $\dim K < n - 2$ or  $0 \in \mathrm{aff}\, K$, (\ref{eq:3.4}) holds by Lemma 5. Therefore, we may assume that $\dim K = n - 2$ and
$0 \notin \mathrm{aff}\, K$. Moreover, by (\ref{eq:1.2}), we may also assume that $K \subseteq e_n^{\bot}$.

Let $T \subseteq e_n^{\bot}$ be an $(n-2)$-dimensional simplex with $0 \notin \mathrm{aff}\, T$.
Choose $\psi\in\mathrm{GL}(n-1)$
with $\det \psi>0$ such that the linear map
$$\phi=\begin{pmatrix}\psi&0\\ 0& 1\end{pmatrix}$$
satisfies $T=\phi S$, where
$S=\mathrm{conv}\{e_1,\ldots,e_{n-1}\}$. A simple calculation shows that $\det \psi = (n-1)!\,
\mathrm{vol}_{n-1}(T_{\mathrm{o}})$. Since $\Phi S \subseteq \mathrm{span}\{e_n\}$ by Lemma 5 (ii), we obtain from (\ref{eq:1.2}) that
\[ \Phi\, T =(\det\phi) \phi^{-t} \Phi S = (n-1)!\, \mathrm{vol}_{n-1}(T_{\mathrm{o}})\, [-a\, e_n,b\,e_n],\]
with constants $a,b\in\mathbb{R}$ depending only on $\Phi$.
Choosing $\vartheta \in \mathrm{SO}(n)$ such that $\vartheta e_n
= -e_n$ and $\vartheta S = S$ (this is possible since $n\geq3$),
it follows from (\ref{eq:1.2}) that \linebreak $\Phi S=\vartheta
\Phi S$. Hence, $a = b \geq 0$ and (\ref{eq:3.4}) holds for $T$.

Now suppose that $P \subseteq e_n^{\bot}$ is a polytope of
dimension $n-2$ such that \linebreak $0 \notin \mathrm{aff}\,P$.
We dissect $P$ into $(n-2)$-dimensional simplices $P
=S_1\cup\ldots\cup S_m$ such that $\dim S_i \cap S_j < n - 2$. As
in the proof of Proposition 6, an application of the
inclusion-exclusion principle and Lemma 5 (i) shows that
\[\Phi P=\Phi S_1+\ldots+\Phi S_m=c\, \mathrm{vol}_{n-1}(P_{\mathrm{o}})\,[-e_n,e_n]\]
for some constant $c \geq 0$ depending only on $\Phi$. Since $\Phi$ is continuous, we conclude that
(\ref{eq:3.4}) holds for all convex bodies $K \subseteq e_n^{\bot}$ of dimension $n - 2$.

\vspace{0.2cm}

Now let $K \in \mathcal{K}^n$ be an arbitrary convex body contained in $u^{\bot}$, $u \in S^{n-1}$.
We want to show that there exist constants $a_1, a_2 \geq 0$ depending only on $\Phi$ such that
\begin{equation}\label{eq:3.5} \Phi K = \left(a_1 \mathrm{vol}_{n-1}(K)+ a_2 \mathrm{vol}_{n-1}(K_{\mathrm{o}}\setminus K)\right)[-u,u].\end{equation}

If $\dim K < n - 1$, then (\ref{eq:3.5}) holds by Lemma 5 and
(\ref{eq:3.4}). Therefore, we may assume that $\dim K = n - 1$.
Moreover, by (\ref{eq:1.2}), we may also assume that $K \subseteq
e_n^{\bot}$.

Let $P \subseteq e_n^{\bot}$ be an $(n-1)$-dimensional polytope. Recall that a face $F$ of $P$ is called visible if
$[0,x)\cap P=\emptyset$ for every $x\in F$. Suppose that $P$
has exactly one visible $(n-2)$-face $F$. Since $\Phi$ is a valuation, we have
\begin{equation} \label{eq117}
\Phi F_{\mathrm{o}}+\Phi P= \Phi P_{\mathrm{o}}+ \Phi F.
\end{equation}
By Theorem 4, there exists a constant $a_1\geq0$ such that $\Phi$
restricts to $a_1 \Pi$ on $\mathcal{K}_{\mathrm{o}}^n$. Hence,
$\Phi F_{\mathrm{o}}=a_1\Pi F_{\mathrm{o}}$ and $\Phi
P_{\mathrm{o}}=a_1\Pi P_{\mathrm{o}}.$  Moreover, from
(\ref{eq:3.4}), it follows that there exists a constant $a_2 \geq
0$ such that
\[\Phi F = a_2 \mathrm{vol}_{n-1}(F_{\mathrm{o}})[-e_n,e_n]=a_2\Pi F_{\mathrm{o}}.\]
Thus, we can rewrite (\ref{eq117}) in the following way
\[a_1\Pi(F_{\mathrm{o}})+\Phi P= a_1\Pi P_{\mathrm{o}}+ a_2 \Pi F_{\mathrm{o}}.\]
Since $\Pi$ is a Minkowski valuation and $\Pi F=\{0\}$, we have
$\Pi P_{\mathrm{o}} = \Pi P + \Pi F_{\mathrm{o}}$. Thus, using
the cancelation law for Minkowski addition, we deduce that
\[\Phi P=a_1 \Pi P +a_2 \Pi F_{\mathrm{o}}=(a_1 \mathrm{vol}_{n-1}(P)+ a_2 \mathrm{vol}_{n-1}(P_{\mathrm{o}}\setminus P))[-e_n,e_n].\]
In order to obtain (\ref{eq:3.5}) for an $(n-1)$-dimensional polytope $P \subseteq e_n^{\bot}$ with $m > 1$ visible $(n-2)$-faces $F_1, \ldots, F_m$, we dissect $P$ as in the proof of Proposition 6.
Finally, since $\Phi$ is continuous, we conclude that
(\ref{eq:3.5}) holds for all convex bodies $K \subseteq e_n^{\bot}$.

In the last step of the proof we combine (\ref{detlow}) and
(\ref{eq:3.5}), to obtain
\begin{equation} \label{fastfertig}
h(\Phi K,u)=a_1 \mathrm{vol}_{n-1}(K|u^\perp)+a_2 \mathrm{vol}_{n-1}((K_{\mathrm{o}}|u^\perp)\setminus (K|u^\perp))
\end{equation}
for every $K$ in $\mathcal{K}^n$ and $u \in S^{n-1}$. Suppose that $a_1 < a_2$. Then, for $u \in S^{n-1}$,
\[h(\Phi K,u)=c_1 \mathrm{vol}_{n-1}(K_{\mathrm{o}}|u^\perp)-c_2 \mathrm{vol}_{n-1}(K|u^\perp),\]
where $c_1, c_2 > 0$. By (\ref{projform}), this is equivalent to
\begin{equation}\label{eq:3.6}h(\Phi K,u)=\int_{S^{n-1}} | u \cdot v |\ d\rho(v), \qquad u\in S^{n-1},\end{equation}
where $\rho$ is the signed Borel measure defined by
\begin{equation} \label{rho17}
\rho = \frac{c_1}{2} S_{n-1}(K_{\mathrm{o}},\cdot)- \frac{c_2}{2}
S_{n-1}(K,\cdot).
\end{equation}
It follows from (\ref{eq:3.6}) that for every polytope $P$ there
exist polytopes $Q_1$ and $Q_2$ such that
\[\Phi P+Q_1=Q_2.\]
Hence, for every polytope $P$, $\Phi P$ is also a polytope.
Since, by (\ref{eq:3.6}), $\Phi P$ is also a generalized zonoid,
we conclude that $\rho\geq 0$, whenever $P$ is a polytope (cf.\
\textbf{\cite[\textnormal{Corollary 3.5.6}]{schneider93}}). But
this is a contradiction, as can be seen by calculating the
surface area measures in (\ref{rho17}) for
$P=\mathrm{conv}\{e_1,\ldots,e_n\}$.

\pagebreak

We conclude that $a_1 \geq a_2$ in (\ref{fastfertig}).
Consequently, there are constants $c_1, c_2 \geq 0$ such that
(\ref{fastfertig}) is equivalent to
\[h(\Phi K,u)=c_1 h(\Pi K,u)+c_2 h(\Pi_{\mathrm{o}} K,u), \qquad u \in S^{n-1},\]
which completes the proof of Theorem 7. \hfill $\blacksquare$

\vspace{1.2cm}

\centerline{\large{\bf{ \setcounter{abschnitt}{5}
\arabic{abschnitt}. A Petty projection type inequality}}}

\vspace{0.7cm} \reseteqn \alpheqn \setcounter{koro}{0}

The fundamental affine isoperimetric inequalities for projection
bodies are the Petty \textbf{\cite{petty67}} projection
inequality and the Zhang \textbf{\cite{zhang91}} projection
inequality: \linebreak Among bodies of given volume polar
projection bodies have maximal volume precisely for ellipsoids
and they have minimal volume precisely for simplices. It is a
major open problem to determine the corresponding results for the
volume of the projection body itself (see, e.g.,
\textbf{\cite{lutwak85, lutwak93}}).

Let $\Phi: \mathcal{K}^n \rightarrow \mathcal{K}^n$ be an
$\mathrm{SL}(n)$ contravariant and homogeneous  continuous
Minkowski valuation. It follows from Theorem 7 that among bodies
$K$ of given volume, the quantity $V(\Phi^*K)$ does not, in
general, attain a minimum. However, it always attains a maximum.
The following result provides a generalization of Petty's
projection inequality. In view of Theorem 7 it is equivalent to
Theorem 3.

\vspace{0.4cm}

\noindent {\bf Theorem 8} \emph{Let $K \in \mathcal{K}^n$ have
nonempty interior. If $\Phi :\mathcal{K}^n \rightarrow
\mathcal{K}^n$ is a nontrivial $\mathrm{SL}(n)$ contravariant and
homogeneous continuous Minkowski \linebreak valuation, then
\begin{equation} \label{petty17}
V(K)^{n-1}V(\Phi^*K) \leq V(B)^{n-1}V(\Phi^*B).
\end{equation}
If $\Phi $ is not a multiple of $\Pi$, there is equality if and
only if $K$ is an ellipsoid containing the origin; otherwise
equality holds if and only if $K$ is an ellipsoid.}

\vspace{0.2cm}

\noindent {\it Proof}\,: Since $\Phi$ is $\mathrm{SL}(n)$ contravariant, we have $\Phi(\vartheta K) = \vartheta \Phi K$ for every $\vartheta \in \mathrm{SO}(n)$ and every $K \in \mathcal{K}^n$.
Thus, there exists an $r_{\Phi} \geq 0$ such that $\Phi B = r_{\Phi} B$. Since inequality (\ref{petty17}) is homogeneous and $\Phi$ is nontrivial, we may assume that $\Phi B = \Pi B=\kappa_{n-1} B$. Here, $\kappa_{m}$ denotes the $m$-dimensional volume of the Euclidean unit ball in $\mathbb{R}^m$.

Since $\Pi_{\mathrm{o}}B=\Pi B$, it follows from Theorem 7 that there exists a constant $0 \leq \lambda \leq 1$ such that for every $K \in \mathcal{K}^n$,
\[\Phi K = \lambda\, \Pi K + (1 - \lambda)\, \Pi_{\mathrm{o}} K. \]
Thus, using (\ref{harsum}), we obtain for convex bodies $K$ with
nonempty interior,
\[\Phi^*K = \lambda \cdot \Pi^*K\, \widehat{+} \, (1- \lambda)\cdot \Pi_{\mathrm{o}}^* K.  \]
An application of the harmonic dual Brunn--Minkowski inequality
(\ref{dualbm}), now yields
\begin{equation} \label{appdualbm}
V(\Phi^*K)^{-1/n} \geq \lambda V(\Pi^*K)^{-1/n} + (1-\lambda)V(\Pi_{\mathrm{o}}^*K)^{-1/n},
\end{equation}
with equality, for $0 < \lambda < 1$, if and only if $\Pi^*K$ and $\Pi_{\mathrm{o}}^*K$ are dilates.

Since $K \subseteq \mathrm{conv}(\{0\} \cup K)$, definition (\ref{defpi}) of the projection operator shows that
$\Pi K \subseteq \Pi_{\mathrm{o}} K$, with equality if and only if $K$ contains the origin. Consequently,
\begin{equation} \label{cont31}
V(\Pi_{\mathrm{o}}^* K)  \leq V(\Pi^*K),
\end{equation}
with equality if and only if $K$ contains the origin. Combining (\ref{appdualbm}) and (\ref{cont31}), yields
\[V(\Phi^*K) \leq V(\Pi^*K).   \]
If $\Phi \neq \Pi$, there is equality if and only if $K$ contains the origin. Finally, an application of the Petty projection inequality (\ref{petty1}), yields (\ref{petty17}) along with the equality conditions.
\hfill $\blacksquare$

\vspace{0.4cm}

The proof of Theorem 8 also shows that the Petty projection inequality is the strongest inequality in the family of inequalities (\ref{petty17}). More precisely:

\vspace{0.4cm}

\noindent {\bf Corollary 9} \emph{Let $K \in \mathcal{K}^n$ have
nonempty interior. If $\Phi :\mathcal{K}^n \rightarrow
\mathcal{K}^n$ is an $\mathrm{SL}(n)$ contravariant and
homogeneous continuous Minkowski valuation such that $\Phi B =
\Pi B$ , then
\begin{equation*}
V(\Phi^*K) \leq V(\Pi^*K).
\end{equation*}
If $\Phi \neq \Pi$, there is equality if and only if $K$ contains
the origin.}

\vspace{0.5cm}

\noindent{\bf Acknowledgments.} The work of the authors was
supported by the \linebreak Austrian Science Fund (FWF), within
the project ``Minkowski valuations and geometric inequalities",
Project number: P\,22388-N13. The second \linebreak author was
also supported by the FWF project ``Analytic and probabilistic
\linebreak methods in Combinatorics", Project number: S9604-N13.

\vspace{0.5cm}

Vienna University of Technology \par Institute of Discrete
Mathematics and Geometry \par Wiedner Hauptstra\ss e 8--10/1046
\par A--1040 Vienna, Austria \\[-0.3cm]
\par franz.schuster@tuwien.ac.at \par thomas.wannerer@tuwien.ac.at

\end{document}